\documentclass[12pt,fleqn]{article}
\usepackage{hyperref}
\usepackage{amsfonts}
\usepackage{amssymb}
\usepackage{graphicx}
\newcommand{\ol}{\setlength{\itemsep}{0pt.}\begin{enumerate}}
\newcommand{\eol}{\end{enumerate}\setlength{\itemsep}{-\parsep}}
\newcommand{\ignore}[1]{}
\newcommand{\quotes}[1]{``#1''}

\newtheorem{THEOREM}{Theorem}[section]
\newenvironment{theorem}{\begin{THEOREM} \hspace{-.85em} {\bf :}
}%
                        {\end{THEOREM}}
\newtheorem{LEMMA}[THEOREM]{Lemma}
\newenvironment{lemma}{\begin{LEMMA} \hspace{-.85em} {\bf :} }%
                      {\end{LEMMA}}
\newtheorem{COROLLARY}[THEOREM]{Corollary}
\newenvironment{corollary}{\begin{COROLLARY} \hspace{-.85em} {\bf
:} }%
                          {\end{COROLLARY}}
\newtheorem{PROPOSITION}[THEOREM]{Proposition}
\newenvironment{proposition}{\begin{PROPOSITION} \hspace{-.85em}
{\bf :} }%
                            {\end{PROPOSITION}}
\newtheorem{DEFINITION}[THEOREM]{Definition}
\newenvironment{definition}{\begin{DEFINITION} \hspace{-.85em} {\bf
:} \rm}%
                            {\end{DEFINITION}}
\newtheorem{EXAMPLE}[THEOREM]{Example}
\newenvironment{example}{\begin{EXAMPLE} \hspace{-.85em} {\bf :}
\rm}%
                            {\end{EXAMPLE}}
\newtheorem{CONJECTURE}[THEOREM]{Conjecture}
\newenvironment{conjecture}{\begin{CONJECTURE} \hspace{-.85em}
{\bf :} \rm}%
                            {\end{CONJECTURE}}
\newtheorem{MAINCONJECTURE}[THEOREM]{Main Conjecture}
\newenvironment{mainconjecture}{\begin{MAINCONJECTURE} \hspace{-.85em}
{\bf :} \rm}%
                            {\end{MAINCONJECTURE}}
\newtheorem{PROBLEM}[THEOREM]{Problem}
\newenvironment{problem}{\begin{PROBLEM} \hspace{-.85em} {\bf :}
\rm}%
                            {\end{PROBLEM}}
\newtheorem{QUESTION}[THEOREM]{Question}
\newenvironment{question}{\begin{QUESTION} \hspace{-.85em} {\bf :}
\rm}%
                            {\end{QUESTION}}
\newtheorem{REMARK}[THEOREM]{Remark}
\newenvironment{remark}{\begin{REMARK} \hspace{-.85em} {\bf :}
\rm}%
                            {\end{REMARK}}
\newtheorem{CLAIM}[THEOREM]{Claim}
\newenvironment{claim}{\begin{CLAIM} \hspace{-.85em} {\bf :}
\rm}%
                            {\end{CLAIM}}
\newcommand{\thm}{\begin{theorem}}
\newcommand{\lem}{\begin{lemma}}
\newcommand{\pro}{\begin{proposition}}
\newcommand{\dfn}{\begin{definition}}
\newcommand{\rem}{\begin{remark}}
\newcommand{\xam}{\begin{example}}
\newcommand{\cnj}{\begin{conjecture}}
\newcommand{\mcnj}{\begin{mainconjecture}}
\newcommand{\prb}{\begin{problem}}
\newcommand{\que}{\begin{question}}
\newcommand{\cor}{\begin{corollary}}
\newcommand{\clm}{\begin{claim}}
\newcommand{\ethm}{\end{theorem}}
\newcommand{\elem}{\end{lemma}}
\newcommand{\epro}{\end{proposition}}
\newcommand{\edfn}{\end{definition}}
\newcommand{\erem}{\bbox\end{remark}}
\newcommand{\exam}{\bbox\end{example}}
\newcommand{\ecnj}{\bbox\end{conjecture}}
\newcommand{\emcnj}{\bbox\end{mainconjecture}}
\newcommand{\eprb}{\bbox\end{problem}}
\newcommand{\eque}{\bbox\end{question}}
\newcommand{\ecor}{\end{corollary}}
\newcommand{\eclm}{\end{claim}}
\newcommand{\beqn}{\begin{equation}}
\newcommand{\eeqn}{\end{equation}}
\newcommand{\bbox}{\begin{flushright} $\Box $ \end{flushright}}
\newcommand{\qed}{\bbox}
\def \1{\mathbf 1}

\def\<{\left<}
\def\>{\right>}
\def \({\left(}
\def \){\right)}

\def \8{\infty}

\newcommand{\remove}[1]{}

\title{The success probability in Lionel Levine's hat problem is strictly decreasing with the number of players, and this is related to interesting questions regarding Hamming powers of Kneser graphs and independent sets in random subgraphs }
\author{Ehud Friedgut \and Gil Kalai\and Guy Kindler}
\begin{document}
\maketitle
\date{}
\begin{abstract}
Lionel Levine's hat challenge has $t$ players, each with a (very large, or infinite) stack of hats on their head, each hat independently colored at random black or white. 
The players are allowed to coordinate before the random colors are chosen, but not after.
Each player sees all hats except for those on her own head.
They then proceed to simultaneously try and each pick a black hat from their respective stacks.
They are proclaimed successful only if they are all correct.
Levine's conjecture was the success probability tends to zero when the number of players grows. 
We prove that this success probability is strictly decreasing in the number of players, and present some connections to questions in graph theory. 

\end{abstract}
\section{Introduction}
The following question proposed by Lionel Levine, arose in the context of his work with Friedrich \cite{FL13}. It gained considerable popularity after being presented in 2011 Tanya Khovanova's blog \cite{Khova}.
Consider $t$ players, each with a stack of $n$ hats on her head, where the hats are chosen independently at random to be either black or white with probability 1/2. Each player sees the hats of every other player, but not her own. Then, simultaneously, all players pick a hat from their respective stacks. The collective of players wins if every single player points to a black hat, else, if even a single player errs, the collective fails. Let $p(t,n)$ be the maximal success probability over all possible strategies that the players can apply. Let $p(t)$ be the limit of $p(t,n)$ as $n$ tends to infinity. The challenge set by Levine was to prove the following conjecture.
\begin{conjecture}\label{Main}
$p(t)$ tends  to 0 as $t$ grows.
\end{conjecture}
In this paper we prove
\begin{theorem}\label{thm:main}
 $p(t+1) < p(t)$ for all $t\geq1$.
 \end{theorem}
 While preparing this paper we were sent a draft of a comprehensive hat-related paper by Buhler, Freiling, Graham, Kariv, Roche, Tiefenbruck, Van Alten and Yeroshkin \cite{BFGKRTVY}, where Theorem \ref{thm:main} is also proven, along with other interesting results and bounds. We refer to their paper as an excellent source for background on the state of the art for this problem. The most prominent landmarks mentioned in their paper are 
$$
0.35 \leq p(2) \leq 0.361607
$$
and
$$
p(t) = \Omega (1/\log(t)),
$$
where the bounds on $p(2)$ are due to them, and the $1/\log(t)$ probably due to Peter Winkler.
The fact that $p(2) \leq 3/8$ was well known folklore in the hatter community, and apparently Alon and Tardos had an approach to improve it \cite{AT}, but not to as tight a bound as 0.362.

 In this paper, we also present some generalizations of Levine's conjecture, relating it to questions regarding independent sets in Hamming products of graphs, and independent sets in random subgraphs.

\section{General setting, strategies and winning sets}
Let us start by defining a general setting that includes the hats game as a special case.
Let $B$ be a fixed ground set, and $\cal{W}$$^{(1)}$ a family of subsets of $B$, that we will call winning sets. (The superscript ${}^{(1)}$ will soon become clear.)
In the corresponding game there are $t \geq 1$ players, each is assigned a point at random from $B$, and each player sees the points the other players were assigned, but not her own point (which is \quotes{on her forehead}). Then, simultaneously, each player choses a winning set, and the collective of players succeeds if every player named a set containing her point.

 Let us now define what a strategy is for $B^t$ (\quotes{the game for $t$ players}), and what a winning set is for such a strategy. For $t=1$ a winning set is any of the sets in  $\cal{W}$$^{(1)}$, and a strategy is a choice of one winning set (or, if you prefer, a function $f: \{\emptyset\} \rightarrow$  $\cal{W}$$^{(1)}$, so the winning set for the strategy $f$ is $f(\emptyset)$.)
 For the $t$-player game a strategy is a $t$-tuple of functions, 
 $(f_1,\ldots,f_t)$, where $f_i : B^{t-1} \rightarrow$ $\cal{W}$$^{(1)}$. 
 For a $t$-tuple $x=(x_1,\ldots,x_t) \in B^t$ let $x^{-i}$ denote the $(t-1)$-tuple obtained from $x$ by deleting the $i$'th coordinate. The winning set for  $(f_1,\ldots,f_t)$ is the set of all $x=(x_1,\ldots,x_t)$ such that for all $i$ it holds that $x_i \in f_i(x^ {-i})$. We denote the family of all winning sets for strategies in the $t$-player game by  $\cal{W}$$^{(t)}$.
 
 An alternative but equivalent way of defining strategies and winning sets for the $t$-player game is the following, which will prove useful for us. (There is a canonical isomorphism between these two definitions). We proceed to define by induction. For $t=1$ we use the previous definition. For $t>1$ we view $B^t$ as  $B^{t-1} \times B$, and let $X_1=B^{t-1}$, and $X_2=B$.  We define a strategy and a winning set for $X_1 \times X_2$.
\\ A strategy for $X_1 \times X_2$ is a pair of functions $f_1,f_2$, with $f_1: X_2 \rightarrow$ $\cal{W}$$^{(t-1)}$, and
$f_2: X_1 \rightarrow$ $\cal{W}$$^{(1)}$.  The winning set for $(f_1,f_2)$ is the set of all points $x_1,x_2$ such that $x_1$ belongs to the winning set $f_1(x_2)$ and $x_2$ belongs to the winning set $f_2(x_1)$. $\cal{W}$$^{(t)}$ is the set of all such winning sets. 

In this paper we will concentrate on $B=\{0,1\}^n$ with various choices for $\cal{W}$$^{(1)}$. We will use the uniform measure on $B^t$, which we denote by $\mu$, and define 
$$
p(t,n) := \max_{W \in {\cal{W}}^t} \mu(W),\ \ \  p(t):= \lim_{n \rightarrow \infty} p(t,n).
$$
The three choices of $\cal{W}$$^{(1)}$ that will interest us are:
\begin{itemize}
\item Let $\cal{W}$$_{dict}$ be the set of dictators, i.e. the set of all $W_i = \{x \in B : x_i=1\}$. This is the basis for the hats game. We will henceforth use $p_{dict}(t,n)$ and $p_{dict}(t)$ for $p(t,n)$ and $p(t)$, the success probabilities in this setting.
 \item Let $\cal{W}$$_{intersecting}$ be the set of all intersecting families in $\{0,1\}^n$, i.e, the set of all $W \subset \{0,1\}^n$ such that if $x,y \in W$ then there exists a coordinate $i$ such that $x_i=y_i=1$. We will use 
 $p_{intersecting}(t,n)$ and $p_{intersecting}(t)$ for the success probabilities in this case.
 \item Let $\cal{W}$ be the set of all balanced monotone families in $\{0,1\}^n$, i.e, all $W$ containing precisely half the points in $\{0,1\}^n$ such that if $x \in W$ and $y_i \geq x_i$ for all $i$ then $y \in W$.  
 We will use 
 $p_{monotone}(t,n)$ and $p_{monotone}(t)$ for the success probabilities in this case.
 \end{itemize}
 Note that every dictatorship is an intersecting family, and every maximal intersecting family is a balanced monotone family, so 
 $$
 p_{monotone}(t,n)\geq  p_{intersecting}(t,n) \geq  p_{dict}(t,n)
 $$
 and
  $$
 p_{monotone}(t)\geq  p_{intersecting}(t) \geq  p_{dict}(t).
 $$
 Thus the following two conjectures are progressively stronger than Conjecture \ref{Main}
 \begin{conjecture}\label{Kneser}
$    p_{intersecting}(t)$ tends  to 0 as $t$ grows.
\end{conjecture}
\begin{conjecture}\label{Monotone}
$    p_{monotone}(t)$ tends  to 0 as $t$ grows.
\end{conjecture}
 \subsection{Winning sets as independent sets in Hamming products of graphs}
Having described the general setting we would like to point out that Conjecture \ref{Kneser} is actually a statement in graph theory. To that end, here are some definitions. 
\dfn
The Kneser graph $K(n)$ is a graph on vertex set $\{0,1\}^n$, with an edge between $x$ and $y$ if $x$ and $y$ have disjoint support, i.e. there is no $i$ for which $x_i=y_i=1$.
\edfn
\dfn
The Hamming product of graphs $G$ and $H$ has vertex set $V(G) \times V(H)$, and an edge between
$(x,v)$ and $(y,u)$ if either $x=y$ and $\{v,u\}$ is an edge in $H$, or $v=u$ and $\{x,y\}$ is an edge in $G$.
We denote it by $G \square H$. There is a canonical isomorphism between $G \square (H \square M)$ and 
$(G \square H) \square M$, so we will treat this product as an associative relation, and write $G^{\square t}$ to denote the $t$-fold Hamming product of $G$ with itself.
\edfn
\dfn
Let $\alpha(G)$ be the size of the largest independent set in $G$, and $\bar{\alpha}(G) := \frac{\alpha(G)}{|V(G)|}$
\edfn
Note that an independent set in $G^{\square t}$ is a subset of $(V(G))^t$ such that its intersection with every 1 dimensional fiber of $(V(G))^t$ is an independent set in $G$. Also note that an independent set in $K(n)$ is an intersecting family. Thus,
\\ {\bf Observation}: 
$$
p_{intersecting}(t,n) = \bar{\alpha}(K(n)^{\square t} ).
$$
So, we may restate Conjecture \ref{Kneser} as
\begin{conjecture}\label{Kneser2}
$$
\lim_{t \rightarrow \infty} \bar{\alpha}(K(n)^{\square t} ) = 0 .
$$
\end{conjecture}
\subsection{Relating the maximal winning set in $B^{t+1}$ to the maximal winning set in a random subset of $B^t$}
 We now return to the general setting of a game on $B^t$ and $B^{t+1}$ and proceed to express $p(t+1)$ as the expected measure of the largest intersection of a winning set and a random set in $B^t$. 
 
 Consider the game on $B^{t+1} = B^t \times B$ and a strategy $(f_1,f_2)$, with $f_1:B \rightarrow {\cal{W}}^{(t)}$ and
  $f_2: B^t \rightarrow {\cal{W}}^{(1)}$. We claim that for a given $f_2$ it is simple to describe an optimal choice of $f_1$.
  Let $ {\cal{W}}^{(1)} = \{W^{(1)}_i\}_{i=1}^r$, and
  first, note that $f_2$ induces a partition of $B^t$ into $V_1,\ldots V_r$, where $V_i:= f_2^{-1}(W^{(1)}_i)$.
  Secondly, note that a random uniform choice of $y \in B$ induces a random subset $R_y \subseteq [r]$ according to the winning sets that $y$ belongs to, i.e.
  $$
  R_y := \{ i: y \in W^{(1)}_i\}.
  $$
  Now, given $f_2$, and a random choice of $x_2 \in B$, how best to define $f_1(x_2)$? Well, observe that a necessary condition for $x_1,x_2$ to be contained in a winning set is for $x_2 \in f_2(x_1)$ which means that $x_2 \in V_i$ for some $i \in I_{x_{2}}$. Therefore, the best choice of $W^{(t)}_j$ for defining $f_1(x_2)=W^{(t)}_j$ is such that it maximizes the probability that a random choice of $x_1 \in B^t$   lands in  $W^{(t)}_j \cap \bigcup_{i \in R_{x_{2}}} V_i$. This implies
 \begin{lemma}\label{p(t+1)}
  $$
 p(t+1)= \max_{B^t=\cup_{i=1}^r V_i } E_{x_2 \in B}\left[\max_{W^{(t)} \in {\cal{W}}^{(t)} } 
 \mu(W \cap \bigcup_{i \in R_{x_{2}}} V_i)\right].
 $$
  \end{lemma}
 Here the first maximum is over all partitions of $B^t$, (each corresponding to a choice of $f_2$), and the second maximum represents the success probability, given $x_2 \in B$, for the optimal choice of $f_1(x_2)$.

 \subsection{A special case: maximal indepedent sets in random subsets of hamming powers of the Kneser graph}
 We can use Lemma \ref{p(t+1)} to find an upper bound for  $p_{intersecting}(t+1,n)$ (and hence for $p_{dict}(t+1,n)$)
 in graph theoretic terms. Let ${\cal{W}}={\cal{W}}_{intersecting}=\{W_i\}_{i=1}^r$ be the family of maximal independent sets in $K(n)$, or, in other words, the family of maximal intersecting sets in $\{0,1\}^n$. A choice of a random vertex $v$ 
 in $K(n)$ induces a choice of a random set $R_v \subseteq [r]$, consisting of all indices $i$ such that $v$ belongs to $W_i$,
 $$
 R_v:=\{i : v \in W_i\}.
 $$
 Each $W \in {\cal{W}}$ has measure 1/2, so the marginal probability of each $i$ belonging to $R_v$ is precisely 1/2.
 Due to a correlation inequality (Kleitman's [K] or Harris's [H]  or FKG [FKG]) these events are non-negatively correlated, i.e., for every $i \in [r], J \subseteq [r]$
 $$
 Pr[i \in R_v | J \subseteq R_v] \geq 1/2.
 $$
 Let ${\cal{D}} = \cup_r {\cal{D}}_r$ denote the set of all such distributions (for all values of $r$).
 We have, then, the following corollary of Lemma \ref{p(t+1)}.
 \beqn\label{1}
 p_{intersecting}(t+1,n)  \leq 
  \max_{r, D \in {\cal{D}}_r, (K(n))^{\square t}=\cup_{i=1}^r V_i} 
  E_{R \sim D}\left[\max_{W \in {\cal{W}}} \mu(W \cap \bigcup_{i \in R} V_i)\right].
 \eeqn
 {\bf Remarks:}
 \begin{itemize}
 \item Equation (\ref{1}) bounds the size of the maximal independent set in the $(t+1)$'th Hamming power of the Kneser graph in terms of the maximal independent set contained in a random subset of the vertices of the $t$'th power. We will expand below on this theme, and raise some conjectures regarding this setting in general graphs.
 \item Recalling that    $p_{intersecting}(t,n) \geq p_{dictator}(t,n)$ makes this approach relevant to solving the hats problem
 \item A similar inequality holds for $ p_{monotone}(t+1,n)$, since monotone increasing subsets, like intersecting fmailies, are positively correlated.
 \end{itemize}

 \subsection{Maximal independent sets in random subgraphs}
 We would like to make a general conjecture regarding independent sets in random subgraphs, that if true, using (\ref{1}), would imply that $ \bar{\alpha}(K(n)^{\square t} )$ tends to 0 as $t$ grows, and thus also prove Levine's conjecture regarding the hats problem, Conjecture \ref{Main}.
 \\ First let us recall some definitions, and make some new ones. 
 \\ For any graph $G$ let $\mu$ denote the uniform measure on $V(G)$.
 \\ Let ${\cal{I}}(G)$ be the family of all independent sets in $G$.
 \\ Let $\bar{\alpha}(G)= \max_{I \in {\cal{I}}(G)} \mu(I)$
 \\ Let ${\cal{D}} = \cup_r {\cal{D}}_r$ denote all distributions on subsets $R$ of some finite set $[r]$, such that every $i \in [r]$ belongs to $R$ independently with probability 1/2, and all these 
 events are positively correlated. 
 \\ Let $\alpha^*(G) = \max_{r, D \in {\cal{D}}_r, V(G)=\bigcup_{i=1}^r V_i}E_{R \sim D} [ \max_{I \in {\cal{I}}(G)} \mu(I \cap (\cup_{i \in R} V_i))]$.
  \\ Let $\alpha^{**}(G) = E_W[ \max_{I \in {\cal{I}}(G)} \mu(I \cap W)]$, where $W$ is chosen uniformly over all subsets of $V(G)$.
  \\ Let $\epsilon^*(\alpha) = \inf_{G: \bar{\alpha}(G)\geq \alpha} \{ \bar{\alpha}(G) - \alpha^*(G)\}$  
  \\ Let $\epsilon^{**}(\alpha) = \inf_{G: \bar{\alpha}(G) \geq \alpha} \{ \bar{\alpha}(G) - \alpha^{**}(G)\}$
 
\begin{conjecture}\label{randomindset}
$\epsilon^*(\alpha)>0$ for all $\alpha \in (0,1/2)$.
 \end{conjecture}
This conjecture would imply Conjecture \ref{Kneser2} (or, equivalently,Conjecture \ref{Kneser})  as follows. Assume, by way of contradiction, that
 $ \bar{\alpha}(K(n)^{\square t} )$ does not tend to 0. Since it is monotone non-increasing in $t$ (this is easy to see, e.g. the success probability of the corresponding game cannot increase with the number of players), and bounded from below, it must tend from above to a limit, say $\alpha$. For large enough $t$ we would have 
$$
  \alpha \leq \bar{\alpha}(K(n)^{\square t} ) < \alpha + \epsilon^*(\alpha)
$$
and hence, using (\ref{1}), and Conjecture \ref{randomindset}

$$
 \bar{\alpha}(K(n)^{\square t+1} ) \leq \alpha^*(K(n)^{\square t}) =
 $$
$$
   =\bar{\alpha}(K(n)^{\square t} ) - \epsilon^*( \bar{\alpha}(K(n)^{\square t} )) <
 \alpha + \epsilon^*(\alpha) - \epsilon^*( \bar{\alpha}(K(n)^{\square t} )) \leq \alpha,
$$
(because $\epsilon(\alpha)$ is non-decreasing. )
So $ \bar{\alpha}(K(n)^{\square t+1} < \alpha$, contradiction.

We do not know that Conjecture \ref{randomindset} is true even in the special cases where the distribution of $I$ is simply binomial (i.e. all $i$ belong to $I$ independently with probability 1/2), or in the even more restricted case where  
each $V_i$ consists of a single vertex. Let us state this last case as a separate conjecture, as it is the purest graph-theoretic statement in this paper, and seems to be of independent interest.
\begin{conjecture}\label{alphastar}
$\epsilon^{**}(\alpha)>0$ for all $\alpha \in (0,1/2)$.
In other words:
\\ There exists a monotone non-decreasing function $\epsilon^{**}: (0,1/2) \rightarrow (0,1/2)$ such that the following holds (Where the point is that $\epsilon^{**}>0$).
If $G$ is a graph on $n$ vertices with maximum independent set of size $\alpha n$, $W$ is a binomial random subset of $V(G)$, and $I_W$ is the maximal independent set contained in $W$, then 
$$
 E_{W}[|I_W|/n] \leq \alpha - \epsilon^{**}(\alpha).
$$
\end{conjecture} 
Noga Alon pointed out to us that one can show, using a random graph, that
$$
\epsilon^{**}(\alpha)< \alpha 2^{-\Omega(\frac{1}{\alpha})}.
$$ 
\section{Blockers and proof of the main theorem}
\subsection{Bounding $p(t+1)$ using blockers}
\remove
{Here is another intriguing conjecture regarding graphs that implies Conjecture \ref{randomindset} and Conjecture \ref{alphastar}. It seems to us almost too good to be true, but we could not find a counterexample.
\begin{conjecture}\label{blockingsets}
For every $\alpha \in (0,1/2)$ there exists  $k$ and $\tau>0$ such that if $G$ is a graph on $n$ vertices with maximal independent set of size $\alpha n$, then there exist disjoint subsets of $V(G)$, $A_1,\ldots A_r$ such that
\begin{enumerate}
\item $|A_i| =k$ for all $i$.
\item $|\bigcup A_i| = \tau n$
\item Every maximal independent set in $G$ intersects every set $A_i$.
\end{enumerate}
\end{conjecture}
The implication from this conjecture to conjecture \ref{randomindset} stems from the simple observation that if $I$ is a random union of positively correlated subsets, then any set $A_i$ of size $k$ has probability at least $2^{-k}$ to be disjoint from $I$, so that we expect many of these sets to be disjoint from $I$, forcing all indepedent sets (simultaneously) to contain many vertices outside $I$. We will work through this argument more precisely below. To that end, let's return to the setting of the hats.
}
We now focus on the hats game, i.e. consider $B =\{0,1\}^n$, with winning sets $W_i = \{x \in B : x_i=1\}$ for $i=1...n$.
Let $\mu$ denote the uniform measure on $B$, and by abuse of notation, also on $B^t$.
 Call a subset of $B^t$ a winning set, if it is the winning set of any strategy for the corresponding game. Write $p(t,n)$ and $p(t)$ for short for $p_{dict}(t,n)$ and $p_{dict}(t)$
\dfn
A blocker $A \subset B^t$ is a set of points that intersects every winning set. 
\edfn

\begin{lemma}\label{blockers}
If there exist disjoint blockers $A_1,\ldots, A_r \subset B^t$, such that
\begin{enumerate}
\item $|A_i| =k$ for all $i$.
\item $\mu(\bigcup A_i) = \beta$
\end{enumerate}
Then $p(t+1) \leq p(t)-2^{-2k-2}\beta/k$.
\end{lemma}
{\bf Proof} : By Lemma \ref{p(t+1)} we know that $p(t+1)$ is bounded by the expectation of the maximal intersection of any winning set in $B^t$ with $V$, a random binomial union of subsets of $B^t$. 
 Now, every blocker  $A$ disjoint from $I$ means that {\em every} winning set contains at least one corresponding point from $A$ (hence missed by $I$), 1 point being $1/k$ of the measure of $A$. If the union of all the blockers missed by $V$ has measure $\tau$ this means every winning set contains a set of measure at least $\tau/k$ that's disjoint from $V$. Now, since every blocker is missed with probability at least $\gamma:=2^{-k}$, the expected proportion of missed blockers is at least $\gamma$ which means that with probability at least $\gamma/2$ a proportion of at least $\gamma/2$ of them is missed, contributing $(\gamma/2)^2\beta$ to the expected measure of the union of missed blockers, which contributes   $(\gamma/2)^2\beta/k$ to the expected measure loss of the maximum that defines $p(t+1)$.
To be concrete: with probability at least $(\gamma/2)$ {\bf every} winning set contains some set of points of measure 
$(\gamma/2)\beta/k$ which falls outside $V$  .

\subsection{Constructing blockers for the hats game}
In this subsection we consider the hats game, and construct, for every $t$, a set of blockers for $B^t$.
This, together with Lemma \ref{blockers}, will prove the main claim of this paper, Theorem \ref{thm:main}
\begin{lemma}
 Let $k(1)=2$, and for $d\ge1$ define $k(d+1):= k(d) {{2k(d)}\choose{k(d)}}$. (So $k$ grows as a tower function of $d$).
 Then, for every $d$ there exist a family of blockers $A_1,\ldots,A_r \subset B^d$ with
 \begin{enumerate}
\item $|A_i| =k(d)$ for all $i$.
\item $\mu(\bigcup A_i) = \frac{2}{k(d)}(1-o(1)).$
\end{enumerate}
 \end{lemma}

\begin{corollary}
For all $d > 1$ 
$$
p(d+1) \le p(d)-\left(\frac{2^{-2k(d)-1}}{k(d)^2}\right)(1-o(1)).
$$
\end{corollary}

{\bf Proof of Lemma}
We will build the family of blockers for $B^d$ inductively.
For $d=1$ the set of blockers for $B$ is the set of all pairs $\{x,\bar{x}\}$. Every dictator must contain precisely one element from each pair. Now, assume we have a family of blockers of size $k(d)$ for $B^d$ as desired.
Let $\ell : = {{2k(d)}\choose{k(d)}}$. We will choose randomly (in a manner to be described below) a series of unordered $\ell$-tuples
$Y^{(j)}=\{y^{(j)}_1,y^{(j)}_2,\ldots,y^{(j)}_\ell\}$, with $y^{(j)}_i \in B$, for $j = 1,2,\ldots$ until the measure of the union of these $\ell$-tuples in $B$ is 
$\frac{1}{\ell}(1-o(1))$. The new blockers for $B^{d+1}$ will be all cartesian products of the form $b \times Y^{(j)}$ where $b$ is one of the blockers we designed for $B^d$. The claims regarding the size of the blockers and the measure of their union are immediate. We must check two things. First, that the product  $Y^{(j)} \times b$ is indeed a blocker for $B^{d+1}$, secondly, that one can choose the desired number of disjoint $\ell$-tuples.
\\ To this end, let us describe how the $\ell$-tuples are formed. We take a random partition of the $n$ coordinates of $B$ into $2d(k)$ sets $S_1,\ldots,S_{2d(k)}$, uniformly over all such partitions. For every $I:= \{ i_1,\ldots,i_{d(k)} \}\subset [2d(k)]$, let $y_I \in B$ be the vector whose 1-support is precisely $S_{i_1}\cup S_{i_2}\cup\ldots\cup S_{i_{d(k)}}$. This defines $\ell$ different vectors corresponding to the specific partition. We proceed to choose such  $\ell$-tuples sequentially at random, and discard any $\ell$-tuple that is not disjoint from all its predecessors. Note that the marginal distribution of every $Y_I$ is uniform, hence if the union of all predecessors of $y^{(j)}$ has measure $\epsilon$ then with probability at least $1-\ell \epsilon$ it will be disjoint from its predecessors, so we may, as claimed, continue until the measure of the union of all $\ell$-tuples is $\frac{1}{\ell}(1-o(1))$.
Finally, we prove that $b \times Y^{(j)} $ is a blocker for $B^d \times B$, where $Y^{(j)}$ is an $\ell$-tuple corresponding to some partition, and $b = \{ x_1,\ldots, x_{k(d)} \}$ is a blocker for $B^d$. 
Let $(f,g)$ be a strategy for $B^d \times B$. For every $x \in B^d$, $g$ picks a dictatorship $g(x)=W_i$,
so define $j_1,\ldots,j_{k(d)}$ by $g(x_i)=W_{j_{i}}$.
These are, respectively, the dictators that the second player guesses when she sees one of the $x$'s from $b$ on the first player's forehead. For $i = 1,\ldots, k(d)$ let $S_{r_i}$ be the part of the partition of $[n]$ (used to define $Y^{(j)}$) that contains $j_i$, and let $I = \{i_{r_1},\ldots,i_{r_{k(d)}}\}$ (or an arbitrary  set of size $k(d)$ containing it if the elements in it are not disjoint).  So the $\ell$-tuple $Y^{(j)}$ contains a vector $y_I$ for which all the coordinates $j_1,\ldots,j_{k(d)}$ are equal to 1, meaning that if the second player has $y_I$ on her forehead she will make a correct guess if the first player has any of the $x$'s in $b$ on her forehead, i.e. 
$$
y_I \in \cap_{i=1}^{k(d)} W_{j_i} =\cap_{x_i \in b} g(x_i)
$$
Now, let $f(y_I)$ be the corresponding winning set of $B^d$ that the first player guesses when seeing $y_I$. By the fact that $b$ is a blocker there exists $x \in b$ that belongs to $f(y_I)$, (and, as mentioned, $y_I$ belongs to the winning set $g(x)$), so the pair $(x,y_I)$ belongs to the winning set of $(f,g)$ - i.e. we have proven that $b\times Y^{(j)}$ intersects every winning set for $B^d \times B$ \qed  
\subsection{Blockers in general graphs?}
In light of the partial success in the previous section, a tempting approach to conjecture \ref{randomindset} is to try and prove the existence of a family of disjoint blockers in any graph $G$, where their size and the measure of their union is a function of $\bar{\alpha}(G)$ . However, this is too good to be true, as pointed out to us by Noga Alon, who came up with the following example. Consider the \quotes{Shift Graph}: Fix a positive integer $m$, the vertices are 
$[m] \times [m]$ , and the edges all pairs of the form $\{(i,j),(j,k)\}$ with $i \neq k$. The maximal independent sets are of the form $A \times B$ where $(A,B)$ is an equi-partition of $[m]$, so they have measure 1/4. But the size of the smallest blocker is $m/2$. It is interesting to ask whether this example is asymptotically the worst case, i.e. is it true that in a graph on $n$ vertices with maximal indepedent sets of size $\Omega(n)$ there are always blockers of size $O(\sqrt{n})$. This strengthens a conjecture of
Bollob\'as, Erd\H{o}s and Tuza, who hypothesized the existence of a blocker of size $o(n)$, (see problem 8 in \cite{E}).

\section{Acknowledgments} We thank Noga Alon and Wojtech Samotij for many useful discussions.

 \end{document}